\documentclass{amsart}

\theoremstyle{definition}

\theoremstyle{remark}

\numberwithin{equation}{section}

\begin{document}

\title{ON THE THEOREM OF M. GOLOMB}

\author{Vugar E. Ismailov}

\address{Mathematics and Mechanics Institute,
Azerbaijan National Academy of Sciences, Az-1141, Baku,
Azerbaijan}

\email{vugaris@mail.ru}

\thanks{This research was partially supported by INTAS under Grant
06-1000015-6283}



\begin{abstract}
Let $X_{1},...,X_{n}$ be compact spaces and $X=X_{1}\times \cdots \times X_{n}.$ Consider the approximation of a function $f\in C(X)$ by sums $g_{1}(x_{1})+\cdots g_{n}(x_{n}),$ where $g_{i}\in C(X_{i}),$ $i=1,...,n.$ In [8], M.Golomb obtained a formula for the error of this approximation in terms of measures constructed on special points of $X$, called "projection cycles". However, his proof had a gap, which was pointed out by Marshall and O'Farrell [15]. But the question if the formula was correct, remained open. The purpose of the paper is to prove that Golomb's formula holds in a stronger form.
\end{abstract}

\maketitle

\section{Introduction}

Let $X_{i},i=1,...,n,$ be compact (Hausdorff) topological spaces. Consider the
approximation to a continuous function $f$ on $X=X_{1}\times
\cdots \times X_{n}$ from the manifold

\begin{equation*}
M=\left\{ \sum_{i=1}^{n}g_{i}(x_{i}):g_{i}\in C(X_{i}),~~i=1,...,n\right\} .
\end{equation*}%
The approximation error is defined as the distance from $f$ to $M$:

\begin{equation*}
E(f)\overset{def}{=}dist(f,M)=\underset{g\in M}{\inf }\left\Vert
f-g\right\Vert _{C(X)}.
\end{equation*}

The well-known duality relation says that

\begin{equation*}
E(f)=\underset{\left\Vert \mu \right\Vert \leq 1}{\underset{\mu \in M^{\bot }%
}{\sup }}\left\vert \int\limits_{X}fd\mu \right\vert ,\eqno(1.1)
\end{equation*}%
where $M^{\bot }$ is the space of regular Borel measures
annihilating all functions in $M$ and $\left\Vert \mu \right\Vert$
stands for the total variation of a measure $\mu $. It should be
also noted that the $\sup $ in
(1.1) is attained by some measure $\mu ^{\ast }$ with the total variation $%
\left\Vert \mu ^{\ast }\right\Vert =1.$ We are interested in the problem: is
it possible to replace in (1.1) the class $M^{\bot }$ by some subclass of it
consisting of measures of simple structure? For the case $n=2,$ this problem
was first considered by Diliberto and Straus [4]. They showed that the
measures induced by so-called "closed lightning bolts" are sufficient for
the equality (1.1).

Let $X=X_{1}\times X_{2}$ and $\pi _{i}$ be the projections of $X$
onto $X_{i},$ $i=1,2.$ A lightning bolt (or, simply, a bolt) is a
finite ordered set $\{a_{1},...,a_{k}\}$ contained in $X$, such
that $a_{i}\neq a_{i+1}$, for $i=1,2,...,k-1$, and either $\pi
_{1}(a_{1})=\pi _{1}(a_{2}),$ $\pi _{2}(a_{2})=\pi _{2}(a_{3})$,
$\pi _{1}(a_{3})=\pi _{1}(a_{4}),...,$ or $\pi _{2}(a_{1})=\pi
_{2}(a_{2}),$ $\pi _{1}(a_{2})=\pi _{1}(a_{3})$, $\pi
_{2}(a_{3})=\pi _{2}(a_{4}),...$ A bolt $\{a_{1},...,a_{k}\}$ is
said to be closed if $k$ is an even number and the set
$\{a_{2},...,a_{k},a_{1}\}$ is also a bolt. These objects have
been exploited in a great deal of works devoted to the uniform
approximation of bivariate functions by univariate functions or the related problems,
though sometimes they appeared under the different names (see,
e.g., [2-7,9-11,13-16,18]). In [4], they were called "permissible lines".
The term "lightning bolt" is due to Arnold [1].

Let $l=\{a_{1},...,a_{2k}\}$ be a closed bolt. Consider a measure $\mu _{l}$
having atoms $\pm \frac{1}{2k}$ with alternating signs at the vertices of $l$. That is,

\begin{equation*}
\mu _{l}=\frac{1}{2k}\sum_{i=1}^{2k}(-1)^{i-1}\delta _{a_{i}}\text{ \ or \ }%
\mu _{l}=\frac{1}{2k}\sum_{i=1}^{2k}(-1)^{i}\delta _{a_{i}},
\end{equation*}%
where $\delta _{a_{i}}$ is a point mass at $a_{i}.$ It is clear that $\mu
_{l}\in M^{\bot }$ and $\left\Vert \mu _{l}\right\Vert \leq 1$. $\left\Vert
\mu _{l}\right\Vert =1$ if and only if the set of vertices of the bolt $l$
having even indices does not intersect with that having odd indices. The
following duality relation was first established by Diliberto and Straus [4]

\begin{equation*}
E(f)=\underset{l\subset X}{\sup }\left\vert \int\limits_{X}fd\mu
_{l}\right\vert ,\eqno(1.2)
\end{equation*}%
where $X=X_{1}\times X_{2}$ and the $\sup $ is taken over all closed bolts
of $X$. In fact, Diliberto and Straus obtained the formula (1.2) for the
case when $X$ is a rectangle in $\mathbb{R}^{2}$ with sides parallel to the
coordinate axis. The same result was independently proved by Smolyak (see
[18]). Yet another proof of (1.2), in the case when $X$ is a Cartesian
product of two compact Hausdorff spaces, was given by Light and Cheney
[14]. For $X$'s other than a rectangle in $\mathbb{R}^{2}$, the theorem
under some additional assumptions appeared in the works [9,11,15]. But we
shall not discuss these works here.

Golomb's paper [8] made a start of a systematic study of approximation of
multivariate functions by various compositions, including sums of univariate
functions. Golomb generalized the notion of a closed bolt to $n$-dimensional
case and obtained the analogue of formula (1.2) for the error of
approximation from the manifold $M$. The objects introduced in [8] were
called "projection cycles" and they are simply sets of the form

\begin{equation*}
p=\{b_{1},...,b_{k};~c_{1},...,c_{k}\}\subset X,\eqno(1.3)
\end{equation*}
with the property that $b_{i}\neq c_{j}$, $i,j=1,...,k$ and for all $\nu
=1,...,n,$ the group of the $\nu $-th coordinates of $c_{1},...,c_{k}$ is a
permutation of that of the $\nu $-th coordinates of $b_{1},...,b_{k}.$ Some
points in the $b$-part $\left( b_{1},...,b_{k}\right) $ or $c$-part $\left(
c_{1},...,c_{k}\right) $ of $p$ may coincide. The measure associated with $p$
is

\begin{equation*}
\mu _{p}=\frac{1}{2k}\left( \sum_{i=1}^{k}\delta
_{b_{i}}-\sum_{i=1}^{k}\delta _{c_{i}}\right).
\end{equation*}

It is clear that $\mu _{p}\in M^{\bot }$ and $\left\Vert \mu _{p}\right\Vert
=1.$ Besides, if $n=2,$ then a projection cycle is the union of closed bolts
after some suitable ordering of its points. Golomb's result states that

\begin{equation*}
E(f)=\underset{p\subset X}{\sup }\left\vert \int\limits_{X}fd\mu
_{p}\right\vert ,\eqno(1.4)
\end{equation*}%
where $X=X_{1}\times \cdots \times X_{n}$ and the $\sup $ is taken
over all projection cycles of $X$. It can be proved that in the
case $n=2,$ the formulas (1.2) and (1.4) are equivalent.
Unfortunately, the proof of (1.4) had a gap, which was many years
later pointed out by Marshall and O'Farrell [15]. But the problem if the formula (1.4) was correct, remained unsolved (see also the more recent monograph by Khavinson [11]).

In Section 2, we will construct families of normalized measures (that is,
measures with the total variation equal to $1$) on projection cycles. Each
measure $\mu _{p}$ defined above will be a member of some family. We will also consider minimal projection cycles and
measures constructed on them. By properties of these measures, we
show that Golomb's formula (1.4) is valid and even in a stronger form.

\bigskip

\section{Measures supported on projection cycles}

First we are going to give an equivalent definition of a projection cycle.
This will be useful in constructing of measures of simple structure and with
the capability to approximate an arbitrary measure in $M^{\bot }$.

In the sequel, $\chi _{a}$ will denote the characteristic function of a
single point set $\{a\}\subset \mathbb{R}$.

\bigskip

\textbf{Definition 2.1.} \textit{Let $X=X_{1}\times \cdots \times
X_{n}$ and $\pi _{i}$ be the projections of $X$ onto the sets $%
X_{i},$ $i=1,...,n.$ We say that a set $p=\{x_{1},...,x_{m}\}\subset X$ is a
projection cycle if there exists a vector $\lambda =(\lambda
_{1},...,\lambda _{m})$ with the nonzero real coordinates such that}

\begin{equation*}
\sum_{j=1}^{m}\lambda _{j}\chi _{\pi _{i}(x_{j})}=0,\text{ \ }i=1,...,n.\eqno%
(2.1)
\end{equation*}

\bigskip

Let us give some explanatory notes concerning Definition 2.1. Fix for a
while the subscript $i.$ Let the set $\{\pi _{i}(x_{j})$, $j=1,...,m\}$ have $%
s_{i}$ different values, which we denote by $\gamma _{1}^{i},\gamma
_{2}^{i},...,\gamma _{s_{i}}^{i}.$ Then (2.1) implies that

\begin{equation*}
\sum_{j}\lambda _{j}=0,
\end{equation*}%
where the sum is taken over all $j$ such that $\pi _{i}(x_{j})=\gamma
_{k}^{i},$ $k=1,...,s_{i}.$ Thus for the fixed $i$, we have $s_{i}$
homogeneous linear equations in $\lambda _{1},...,\lambda _{m}.$ The
coefficients of these equations are the integers $0$ and $1.$ By varying $i$%
, we obtain $s=\sum_{i=1}^{n}s_{i}$ such equations. Hence (2.1), in its
expanded form, stands for the system of these equations. One can observe
that if this system has a solution $(\lambda _{1},...,\lambda _{m})$ with
the nonzero real components $\lambda _{i},$ then it also has a solution $%
(n_{1},...,n_{m})$ with the nonzero integer components $n_{i},$
$i=1,...,m.$ This means that in Definition 2.1, we can replace the
vector $\lambda $ by the vector $n=(n_{1},...,n_{m})\,$, where
$n_{i}\in \mathbb{Z}\backslash \{0\},$ $i=1,...,m.$ Thus, Definition 2.1 is equivalent to the following definition.

\bigskip

\textbf{Definition 2.2.} \textit{A set
$p=\{x_{1},...,x_{m}\}\subset X$ is called a projection cycle if there exist nonzero
integers $n_{1},...,n_{m}$ such that}

\begin{equation*}
\sum_{j=1}^{m}n_{j}\chi _{\pi _{i}(x_{j})}=0,\text{ \ }i=1,...,n.\eqno(2.2)
\end{equation*}

\bigskip

\textbf{Proposition 2.3.} \textit{Definition 2.2 is equivalent to Golomb's definition of a projection cycle.}

\begin{proof} Let $p=\{x_{1},...,x_{m}\}$ be a projection cycle with respect to Definition 2.2.
By $b$ and $c$ denote the set of all points $x_{i}$ such that the integers $%
n_{i}$ associated with them in (2.2) are positive and negative
correspondingly. Write out each point $x_{i}$ $n_{i}$ times if $n_{i}>0$ and
$-n_{i}$ times if $n_{i}<0.$ Then the set $\{b;c\}$ is a projection cycle
with respect to Golomb's definition (see Introduction). The inverse is also
true. Let a set $p_{1}=\{b_{1},...,b_{k};~c_{1},...,c_{k}\}$ be a projection
cycle with respect to Golomb's definition. Here, some points $b_{i}$ or $%
c_{i}$ may be repeated. Let $p=\{x_{1},...,x_{m}\}$ stand for the set $p_{1}$%
, but with no repetition of its points. Let $n_{i}$ show how many
times $x_{i}$ appear in $p_{1}.$ We take $n_{i}$ positive if
$x_{i}$ appears in
the $b$-part of $p_{1}$ and negative if it appears in the $c$-part of $%
p_{1}. $ Clearly, the set $\{x_{1},...,x_{m}\}$ is a projection cycle with
respect to Definition 2.2, since the integers $n_{i},$ $i=1,...,m,$ satisfy
Eq. (2.2).
\end{proof}

In the sequel, we will use Definition 2.1. A pair $\left\langle p,\lambda
\right\rangle ,$ where $p$ is a projection cycle in $X$ and $\lambda $ is a
vector associated with $p$ by (2.1), will be called a "projection
cycle-vector pair" of $X.$ To each such pair $\left\langle p,\lambda
\right\rangle $ with $p=\{x_{1},...,x_{m}\}$ and $\lambda =(\lambda
_{1},...,\lambda _{m})$, we correspond the measure

\begin{equation*}
\mu _{p,\lambda }=\frac{1}{\sum_{j=1}^{m}\left\vert \lambda _{j}\right\vert }%
\sum_{j=1}^{m}\lambda _{j}\delta _{x_{j}}.\eqno(2.3)
\end{equation*}

Clearly, $\mu _{p,\lambda }\in M^{\bot }$ and $\left\Vert \mu _{p,\lambda
}\right\Vert =1$. We will also deal with measures supported on some certain subsets of projection cycles called minimal projection cycles. A projection cycle is said to be minimal if it does not contain any projection
cycle as its proper subset. For example, the set $p=%
\{(0,0,0),~(0,0,1),~(0,1,0),~(1,0,0),~(1,1,1)\}$ is a minimal projection
cycle in $\mathbb{R}^{3},$ since the vector $\lambda =(2,-1,-1,-1,1)$
satisfies Eq. (2.1) and there is no such vector for any other subset of $p$.
Adding one point $(0,1,1)$ from the right to $p$, we will also have a
projection cycle, but not minimal. Note that in this case, $\lambda $ can be
taken as $(3,-1,-1,-2,2,-1).$

\bigskip

\textbf{Remark 1.} A minimal projection cycle under the name of a "loop" was introduced in
the work of Klopotowski, Nadkarni, Rao [12].

\bigskip

To prove our main result we need some auxiliary facts.

\bigskip

\textbf{Lemma 2.4.} (1)\textit{\ The vector $\lambda =(\lambda
_{1},...,\lambda _{m})$ associated with a minimal projection cycle $%
p=(x_{1},...,x_{m})$ is unique up to multiplication by a constant.}

(2)\textit{\ If in (1), $\sum_{j=1}^{m}\left\vert \lambda _{j}\right\vert =1,
$ then all the numbers $\lambda _{j}$, $j=1,...,m,$ are rational.}

\begin{proof} Let $\lambda ^{1}=(\lambda _{1}^{1},...,\lambda
_{m}^{1})$ and $\lambda ^{2}=(\lambda _{1}^{2},...,\lambda
_{m}^{2})$ be any two vectors associated with $p.$ That is,

\begin{equation*}
\sum_{j=1}^{m}\lambda _{j}^{1}\chi _{\pi _{i}(x_{j})}=0\text{ and }%
\sum_{j=1}^{m}\lambda _{j}^{2}\chi _{\pi _{i}(x_{j})}=0,\text{ \ }i=1,...,n.
\end{equation*}%
After multiplying the second equality by $c=\frac{\lambda _{1}^{1}}{\lambda
_{1}^{2}}$ and subtracting from the first, we obtain that

\begin{equation*}
\sum_{j=2}^{m}(\lambda _{j}^{1}-c\lambda _{j}^{2})\chi _{\pi _{i}(x_{j})}=0%
\text{, \ }i=1,...,n.
\end{equation*}%
Now since the cycle $p$ is minimal, $\lambda _{j}^{1}=c\lambda _{j}^{2},$
for all $j=1,...,m.$

The second part of the proposition is a consequence of the first
part. Indeed, let $n=(n_{1},...,n_{m})$ be a vector with the
nonzero integer coordinates associated with $p.$ Then the vector
$\lambda ^{^{\prime
}}=(\lambda _{1}^{^{\prime }},...,\lambda _{m}^{^{\prime }}),$ where $%
\lambda _{j}^{^{\prime }}=\frac{n_{j}}{\sum_{j=1}^{m}\left\vert
n_{j}\right\vert },$ $j=1,...,m,$ is also associated with $p.$ All
the coordinates of $\lambda ^{^{\prime }}$ are rational and
therefore by the first part of the proposition, it is the unique
vector satisfying $\sum_{j=1}^{m}\left\vert \lambda _{j}^{^{\prime
}}\right\vert =1. $
\end{proof}

By this proposition, a minimal projection cycle $p$ uniquely (up to a sign) defines the
measure

\begin{equation*}
~\mu _{p}=\sum_{j=1}^{m}\lambda _{j}\delta _{x_{j}},\text{ \ }%
\sum_{j=1}^{m}\left\vert \lambda _{j}\right\vert =1.
\end{equation*}

\bigskip

\textbf{Lemma 2.5} (see [17]). \textit{Let $\mu $ be a normalized orthogonal measure on a
projection cycle $l\subset X$. Then it is a convex combination of normalized orthogonal measures on minimal
projection cycles of $l$. That is,}

\begin{equation*}
\mu =\sum_{i=1}^{s}t_{i}\mu _{l_{i}},\text{ }\sum_{i=1}^{s}t_{i}=1,~t_{i}>0,
\end{equation*}
\textit{where $l_{i},$ $i=1,...,s,$ are minimal projection cycles in $l.$}

\bigskip

This lemma follows from the result of Navada (see Theorem 2 of [17]): Let $S\subset X_{1}\times
\cdots \times X_{n}$ be a finite set. Then any extreme point of the convex set of measures $\mu $ on $S$, $\mu \in M^{\bot }$, $\left\Vert \mu \right\Vert \leq 1$, has its support on a minimal projection cycle contained in $S$.

\bigskip

\textbf{Remark 2.} In the case $n=2$, Lemma 2.5 was proved by Medvedev (see [11, p.77]).

\bigskip

\textbf{Lemma 2.6} (see [11, p.73]). \textit{Let $X=X_{1}\times
\cdots \times X_{n}$ and $\pi _{i}$ be the projections of $X$ onto
the sets $X_{i},$ $i=1,...,n.$ In order that a measure $\mu \in
C(X)^{\ast }$ be orthogonal to the subspace $M$, it is necessary
and sufficient that}

\begin{equation*}
\mu \circ \pi _{i}^{-1}=0,\text{ }i=1,...,n.
\end{equation*}

\bigskip

\textbf{Lemma 2.7} (see [11, p.75]). \textit{Let $\mu \in M^{\bot }$ and $%
\left\Vert \mu \right\Vert =1.$ Then there exist a net of measures $\{\mu
_{\alpha }\}\subset M^{\bot }$ weak$^{\text{*}}$ converging in $C(X)^{\ast }$
to $\mu $ and satisfying the following properties:}

1) $\left\Vert \mu _{\alpha }\right\Vert =1;$

2) \textit{The closed support of each $\mu _{\alpha }$ is a finite set.}

\bigskip

\textbf{Theorem 2.8.} \textit{The error of approximation from the manifold $%
M $ obeys the equality}
\begin{equation*}
E(f)=\underset{l\subset X}{\sup }\left\vert \int\limits_{X}fd\mu
_{l}\right\vert ,
\end{equation*}
\textit{where the $\sup $ is taken over all minimal projection cycles of
$X.$}

\begin{proof} Let $\overset{\sim }{\mu }$ be a
measure with a finite support $\{x_{1},...,x_{m}\}$ and orthogonal
to the space $M.$
Put $\lambda _{j}=\overset{\sim }{\mu }(x_{j}),$ $j=1,...m.$ By Lemma 2.6, $%
\overset{\sim }{\mu }(\pi _{i}^{-1}(\pi _{i}(x_{j})))=0,$ for all $%
i=1,...,n, $ $j=1,...,m.$ Fix the indices $i$ and $j.$ Then we have the
equation $\sum_{k}\lambda _{k}=0,$ where the sum is taken over all indices $%
k $ such that $\pi _{i}(x_{k})=\pi _{i}(x_{j}).$ Varying $%
i$ and $j,$ we obtain a system of such equations, which concisely can be
written as

\begin{equation*}
\sum_{k=1}^{m}\lambda _{k}\chi _{\pi _{i}(x_{k})}=0,\text{ \ }i=1,...,n.
\end{equation*}%
This means that the finite support of $\overset{\sim }{\mu }$ forms a
projection cycle. Therefore, a net of measures approximating the given
measure $\mu $ in Lemma 2.7 are all of the form (2.3).

Let now $\mu _{p,\lambda }$ be any measure of the form (2.3). Since $\mu
_{p,\lambda }\in M^{\bot }$ and $\left\Vert \mu _{p,\lambda }\right\Vert =1,$
we can write

\begin{equation*}
\left\vert \int\limits_{X}fd\mu _{p,\lambda }\right\vert =\left\vert
\int\limits_{X}(f-g)d\mu _{p,\lambda }\right\vert \leq \left\Vert
f-g\right\Vert ,\eqno(2.4)
\end{equation*}%
where $g$ is an arbitrary function in $M$. It follows from (2.4) that

\begin{equation*}
\underset{\left\langle p,\lambda \right\rangle }{\sup }\left\vert
\int\limits_{X}fd\mu _{p,\lambda }\right\vert \leq E(f),\eqno(2.5)
\end{equation*}%
where the $\sup $ is taken over all projection cycle-vector pairs of $X.$

Consider the general duality relation (1.1). Let $\mu _{0}$ be a measure
reaching the supremum in (1.1) and $\left\{ \mu _{p,\lambda }\right\} $ be a
net of measures of the form (2.3) approximating $\mu _{0}$ in the weak$^{%
\text{*}}$ topology of $C(X)^{\ast }.$ We have already known that this is
possible. For any $\varepsilon >0,$ there exists a measure $\mu
_{p_{0},\lambda _{0}}$ in $\left\{ \mu _{p,\lambda }\right\} $ such that

\begin{equation*}
\left\vert \int\limits_{X}fd\mu _{0}-\int\limits_{X}fd\mu _{p_{0},\lambda
_{0}}\right\vert <\varepsilon .
\end{equation*}%
From the last inequality we obtain that

\begin{equation*}
\left\vert \int\limits_{X}fd\mu _{p_{0},\lambda _{0}}\right\vert >\left\vert
\int\limits_{X}fd\mu _{0}\right\vert -\varepsilon =E(f)-\varepsilon .
\end{equation*}%
Hence,

\begin{equation*}
\underset{\left\langle p,\lambda \right\rangle }{\sup }\left\vert
\int\limits_{X}fd\mu _{p,\lambda }\right\vert \geq E(f).\eqno(2.6)
\end{equation*}%
From (2.5) and (2.6) it follows that

\begin{equation*}
\underset{\left\langle p,\lambda \right\rangle }{\sup }\left\vert
\int\limits_{X}fd\mu _{p,\lambda }\right\vert =E(f).\eqno(2.7)
\end{equation*}

By Lemma 2.5,

\begin{equation*}
\mu _{p,\lambda }=\sum_{i=1}^{s}t_{i}\mu _{l_{i}},
\end{equation*}%
where $l_{i}$, $i=1,...,s,$ are minimal projection cycles in $p$ and $%
\sum_{i=1}^{s}t_{i}=1,~t_{i}>0.$ Let $k$ be an index in the set $%
\{1,...,s\} $ such that

\begin{equation*}
\left\vert \int\limits_{X}fd\mu _{l_{k}}\right\vert =\max \left\{ \left\vert
\int\limits_{X}fd\mu _{l_{i}}\right\vert ,\text{ }i=1,...,s\right\} .
\end{equation*}%
Then

\begin{equation*}
\left\vert \int\limits_{X}fd\mu _{p,\lambda }\right\vert \leq \left\vert
\int\limits_{X}fd\mu _{l_{k}}\right\vert .\eqno(2.8)
\end{equation*}%
Now since

\begin{equation*}
\left\vert \int\limits_{X}fd\mu _{l}\right\vert \leq E(f),
\end{equation*}%
for any minimal cycle $l,$ from (2.7) and (2.8) we obtain the assertion of the
theorem.
\end{proof}

\textbf{Remark 3.} Theorem 2.8 not only proves Golomb's formula, but also improves it. Indeed, based on Proposition 2.3, one can easily observe that
the formula (1.4) is equivalent to the formula

\begin{equation*}
E(f)=\underset{\left\langle p,\lambda \right\rangle }{\sup }\left\vert
\int\limits_{X}fd\mu _{p,\lambda }\right\vert ,
\end{equation*}%
where the $\sup $ is taken over all projection cycle-vector pairs $%
\left\langle p,\lambda \right\rangle $ of $X$ provided that all the numbers $%
\lambda _{i}/ \sum_{j=1}^{m}\left\vert \lambda _{j}\right\vert $, $%
i=1,...,m,$ are rational. But by Lemma 2.4, minimal projection cycles enjoy this property.


\bigskip



\begin{thebibliography}{99}
\bibitem{1} V.I.Arnold, On functions of three variables, \textit{Dokl. Akad.
Nauk SSSR} \textbf{114} (1957), 679-681; English transl, \textit{Amer. Math.
Soc. Transl.} \textbf{28} (1963), 51-54.

\bibitem{2} M-B.A.Babaev, Estimates and ways for determining the exact value
of the best approximation of functions of several variables by
superpositions of functions of a smaller number of variables (Russian),
\textit{Special questions in the theory of functions} \textit{(Russian)},
Izdat. "Elm", Baku, 1977, 3--23.

\bibitem{3} D. Braess and A. Pinkus, Interpolation by ridge functions,
\textit{J.Approx. Theory} \textbf{73} (1993), 218-236.

\bibitem{4} S.P.Diliberto and E.G.Straus , On the approximation of a
function of several variables by the sum of functions of fewer variables,
\textit{Pacific J. Math.} \textbf{1} (1951), 195-210.

\bibitem{5} N.Dyn, W.A.Light and E.W. Cheney, Interpolation by
piecewise-linear radial basis functions, \textit{J. Approx. Theory.} \textbf{%
59} (1989), 202-223.

\bibitem{6} A.L.Garkavi, V.A.Medvedev, S.Ya.Khavinson, On the existence of a
best uniform approximation of functions of two variables by sums of the type
$\varphi (x)+\psi (y)$ , \textit{Sibirskii Mat. Zh.}, \textbf{36 }%
(1995),819-827; English transl. in \textit{Siberian Math. J.}, \textbf{36}
(1995), 707-713.

\bibitem{7} M.v.Golitschek and W.A.Light, Approximation by solutions of the
planar wave equation, \textit{Siam J.Numer. Anal.} \textbf{29} (1992),
816-830.

\bibitem{8} M.Golomb, Approximation by functions of fewer variables, \textit{%
On numerical approximation. Proceedings of a Symposium. Madison 1959. Edited
by R.E.Langer. The University of Wisconsin Press,} 275-327.

\bibitem{9} S.Ja.Havinson, A Chebyshev theorem for the approximation of a
function of two variables by sums of the type $\varphi \left( {x}\right)
+\psi \left( {y}\right) ,$ \textit{Izv. Acad. Nauk. SSSR Ser. Mat.} \textbf{%
33} (1969), 650-666; English tarnsl. \textit{Math. USSR Izv.} \textbf{3}
(1969), 617-632.

\bibitem{10} V.E.Ismailov, Methods for computing the least deviation from the
sums of functions of one variable, \textit{Sibirski Matematicheski Zhurnal}
\textbf{47} (2006), 1076--1082; English transl. in \textit{Siberian
Mathematical Journal} \textbf{47} (2006), 883-888.

\bibitem{11} S.Ya.Khavinson, \textit{Best approximation by linear
superpositions (approximate nomography),} Translated from the Russian
manuscript by D. Khavinson. Translations of Mathematical Monographs, 159.
American Mathematical Society, Providence, RI, 1997, 175 pp.

\bibitem{12} A.Klopotowski, M.G.Nadkarni, K.P.S.Bhaskara Rao, When is $%
f(x_{1},x_{2},...,x_{n})=u_{1}(x_{1})+u_{2}(x_{2})+\cdots +u_{n}(x_{n})$ ?,
\textit{Proc. Indian Acad. Sci. Math. Sci.} \textbf{113 }(2003), 77--86.

\bibitem{13} A.Klopotowski, M.G.Nadkarni, Shift invariant measures and
simple spectrum, \textit{Colloq. Math.} \textbf{84/85} (2000), 385-394.

\bibitem{14} W.A. Light and E.W. Cheney , On the approximation of a
bivariate function by the sum of univariate functions, \textit{J.Approx.
Theory} \textbf{29} (1980), 305-323.

\bibitem{15} D.E.Marshall and A.G.O'Farrell, Approximation by a sum of two
algebras. The lightning bolt principle, \textit{J. Funct. Anal.} \textbf{52}
(1983), 353-368.

\bibitem{16} D.E.Marshall and A.G.O'Farrell. Uniform approximation by real
functions, \textit{Fund. Math.} \textbf{104} (1979), 203-211.

\bibitem{17} K.G.Navada, Some remarks on good sets,
\textit{Proc. Indian Acad. Sci. Math. Sci.} \textbf{114} (2003), No.4, 389--397.

\bibitem{18} Ju.P.Ofman, Best approximation of functions of two variables by
functions of the form $\varphi (x)+\psi (y)$, \textit{Izv. Akad. Nauk. SSSR
Ser.Mat.} \textbf{25} (1961), 239-252; English transl. \textit{Amer. Math.
Soc. Transl.} \textbf{44} (1965), 12-28.
\end{thebibliography}
\end{document}